\documentclass{article}
\usepackage[english]{babel}
\usepackage{layout}
\usepackage[latin1]{inputenc}
\usepackage{mathrsfs}
\usepackage{amsmath}
\usepackage{amsthm}
\usepackage{amssymb}
\usepackage{indentfirst}
\usepackage{ifthen}
\newcommand{\cit}[1]{{\rm \textbf{#1}}}

\newcommand{\Ref}[2]{\cit{%
\ifthenelse{\equal{#1}{thm}}{Theorem}{}%
\ifthenelse{\equal{#1}{prop}}{Proposition}{}%
\ifthenelse{\equal{#1}{lem}}{Lemma}{}%
\ifthenelse{\equal{#1}{cor}}{Corollary}{}%
\ifthenelse{\equal{#1}{defn}}{Definition}{}%
\ifthenelse{\equal{#1}{oss}}{Remark}{}%
\ifthenelse{\equal{#1}{sec}}{Section}{}%
\ifthenelse{\equal{#1}{ex}}{Example}{}%
\ifthenelse{\equal{#1}{ssec}}{Subsection}{}%
\ifthenelse{\equal{#1}{conj}}{Conjecture}{}%
\  \ref{#1:#2}%
}}

\newcommand{\hk}{Hyperk\"{a}hler }

\newcommand{\ktipo}{$K3^{[2]}$-type }
\newcommand{\ie}{i.~e.~}

\theoremstyle{plain} 
\newtheorem{prop}{Proposition}[section]
\newtheorem{thm}[prop]{Theorem}
\newtheorem{lem}[prop]{Lemma} 
\newtheorem{cor}[prop]{Corollary}
\newtheorem{conj}[prop]{Conjecture}
\theoremstyle{remark}
\newtheorem{oss}[prop]{Remark}
\newtheorem{ex}[prop]{Example}
\newtheorem*{dim1}{Proof of \Ref{thm}{main_intro}}
\newtheorem*{dim2}{Proof of \Ref{thm}{fixedpoints}}

\theoremstyle{definition}
\newtheorem{defn}[prop]{Definition}
\title{Symplectic involutions on deformations of $K3^{[2]}$}
\author{Giovanni Mongardi}
\begin{document}

\maketitle

\abstract{Let X be a Hyperk\"{a}hler manifold deformation equivalent to the Hilbert square of a K3 surface and let $\varphi$ be an involution preserving the symplectic form. We prove that the fixed locus of $\varphi$ consists of 28 isolated points and 1 K3 surface, moreover the anti-invariant lattice of the induced involution on $H^2(X,\mathbb{Z})$ is isomorphic to $E_8(-2)$. Finally we prove that any couple consisting of one such manifold and a symplectic involution on it can be deformed into a couple consisting of the Hilbert square of a K3 surface and the involution induced by a symplectic involution on the K3 surface.}

\section{Introduction}

The aim of the present paper is to prove that involutions preserving the symplectic form on a Hyperk\"{a}hler variety of $K3^{[2]}$-type are all deformation equivalent.\\ Many papers on automorphisms of \hk manifolds have appeared in recent years, starting from the foundational work of Nikulin \cite{nik1} and Mukai \cite{muk} and an explicit example of Morrison \cite{mor} in the case of K3 surfaces. Then isolated examples of such automorphisms were given by Namikawa \cite{nami}, by Beauville \cite{beau2} and later by Kawatani \cite{kaw} and Amerik \cite{amerik}. Some further work was done in the case of generalized Kummer varieties by Boissi\`{e}re, Nieper-Wi\ss kirchen and Sarti \cite{boiniesar}. Some general work on automorphisms and birational maps was done by Oguiso \cite{ogu3}, Boissi\`{e}re \cite{boi} and Boissi\`{e}re and Sarti \cite{bs}, while order 2 automorphisms were analyzed by Beauville \cite{beau} and Camere \cite{cam}. Before those works on involutions came the work of O'Grady (\cite{ogr} and \cite{ogr2}) on Double-EPW sextics which are naturally endowed with an antisymplectic involution.\\
In the following $X$ will always be a Hyperk\"{a}hler manifold deformation equivalent to the Hilbert square of a K3 surface (\ie a manifold of $K3^{[2]}$-type), and $\sigma$ will denote a holomorphic symplectic form on $X$. We recall that the integral second cohomology of a \hk manifold is endowed with the Beauville-Bogomolov integral quadratic form. \Ref{sec}{prel} contains preliminaries on \hk manifolds and quadratic forms.\\
Given a bimeromorphic map $\varphi$ of $X$ we remark that $\varphi^*$ is defined on $H^2(X,\mathbb{Z})$ since the indeterminacy locus has codimension at least 2. 
Let $\varphi\,:\,X\,\dashrightarrow\,X$ be a meromorphic involution, \ie $\varphi\circ\varphi=Id$. Then $\varphi^*\sigma=\pm\sigma$. 
In this paper we are interested mainly in symplectic involutions, \ie $\,$involutions $\varphi$ such that
\begin{equation}
\varphi^*(\sigma)=\sigma.
\end{equation}
In \Ref{sec}{gen_lattices} we will prove general results concerning finite automorphism groups of manifolds of \ktipo. An interesting question is whether these groups are induced by finite automorphism groups on $K3$'s. We formalize this question as follows. 
Let $G$ be a group acting faithfully on $X$, we call $(X,G)$ a couple. Two couples $(X,G)$ and $(Y,G)$ are isomorphic if there exists a $G$-equivariant isomorphism $X\rightarrow Y$.  There is a natural notion of deformations of the couple $(X,G)$:
\begin{defn}
Given a manifold $X$ and a group $G$ acting faithfully on it we call a $G$-deformation of $X$ (or a deformation of the couple $(X,G)$) the following data:
\begin{itemize}
\item A flat family $\mathcal{X}\rightarrow B$ and a faithful action of the group $G$ on $\mathcal{X}$ inducing fibrewise faithful actions of $G$. 
\item  A map $\{0\}\rightarrow B$ and a $G$-equivariant isomorphism $\mathcal{X}_0\rightarrow X$.
\end{itemize}
\end{defn}
Let $G$ be a cyclic group generated by the automorphism $\varphi$ of $X$. We will denote by $(X,\varphi)$ the couple $(X,G)$. Let $\psi$ be an automorphism of a $K3$ surface $S$, we denote $\psi^{[2]}$ the automorphism it induces on $S^{[2]}$.
\begin{defn}
Let $X$ be a Hyperk\"{a}hler manifold of $K3^{[2]}$-type endowed with an automorphism of finite order $\psi$. The couple $(X,\psi)$ is \emph{standard} if there exists a K3 surface $S$ and an automorphism $\psi'$ of $S$ such that the couples $(X,\psi)$ and $(S^{[2]},\psi'^{[2]})$ are deformation equivalent.
The couple $(X,\psi)$ is \emph{exotic} if it is not standard.
\end{defn} 
We remark that not all automorphisms of a manifold of $K3^{[2]}$-type are standard. In fact there is an example of Namikawa of an exotic automorphism of order 3 of the Fano variety of lines on a particular cubic fourfold (see \cite{nami}).\\
The main result of the present paper is the following:
\begin{thm}\label{thm:main_intro}
Let $X$ be a manifold of $K3^{[2]}$-type, and let $\varphi\in\,Aut(X)$ be a symplectic involution. Then there exists a K3 surface S endowed with a symplectic involution $\psi$ such that $(X,\varphi)$ and $(S^{[2]},\psi^{[2]})$ are deformation equivalent.
\end{thm}
The above result proves the following conjecture made by Camere:
\begin{conj}\label{conj:camere_c}\cite{cam}
Let $X$ be a \hk manifold deformation equivalent to the Hilbert square of a $K3$ surface, and let $\varphi$ be an involution of $X$ preserving the holomorphic symplectic form. Then the fixed locus $X^{\varphi}$ does not contain complex tori. 
\end{conj}

\section{Preliminaries}\label{sec:prel}
\subsection{Hyperk\"{a}hler manifolds}\label{ssec:pre_hk}
This subsection summarizes a few facts about Hyperk\"{a}hler manifolds, the interested reader can consult \cite{huy}.
\begin{defn}
Let $X$ be a K\"{a}hler manifold, it is called a Hyperk\"{a}hler manifold if the following hold:
\begin{itemize}
\item $X$ is compact.
\item $X$ is simply connected.
\item $H^{2,0}(X)=\mathbb{C}\sigma_X$, where $\sigma_X$ is a symplectic form, that is a holomorphic 2-form which is closed and everywhere nondegenerate.
\end{itemize}
\end{defn}
%We remind that a \hk manifold of dimension 2 is a K3 surface.\\
We remark that the isometry class of $H^2(X,\mathbb{Z})$ with the Beauville-Bogomolov form is invariant under smooth deformations.
\begin{ex}
Let $S$ be a $K3$ surface, then $S^{[2]}$ is a \hk manifold. Furthermore $H^2(S^{[2]},\mathbb{Z})$ endowed with its Beauville-Bogomolov pairing is isomorphic to the lattice
\begin{equation}\label{latticeL}
L=U\oplus U\oplus U\oplus E_8(-1)\oplus E_8(-1)\oplus (-2).
\end{equation}
Here $U$ is the hyperbolic lattice, $E_8(-1)$ is the unique unimodular even negative definite lattice of rank 8, $(-2)$ is $(\mathbb{Z},q)$ with $q(1)=-2$ and $\oplus$ denotes orthogonal direct sum. 
We remark that there exists a class $v\in H^{1,1}(S^{[2]},\mathbb{Z})$ of square $-2$ such that $(v,H^{2}(S^{[2]},\mathbb{Z}))=2\mathbb{Z}$.
\end{ex}
Let X be a Hyperk\"{a}hler manifold of $K3^{[2]}$-type, thus $H^2(X,\mathbb{Z})\,\cong\,L$. A marking $f$ of $X$ is an isometry \\ $f\,:\,H^2(X,\mathbb{Z})\,\rightarrow\,L$. 
Marked Hyperk\"{a}hler manifolds are particularly interesting because there exists a moduli space of marked Hyperk\"{a}hler manifolds and there is a good notion of period map (see \cite[Chapter 25]{huy}). Let $\mathcal{M}_{K3^{[2]}}$ be the moduli space of marked Hyperk\"{a}hler manifolds $(X,f)$ deformation equivalent to $K3^{[2]}$, let
\begin{equation}\nonumber \Omega=\{\omega\,\in\,\mathbb{P}(L\otimes\mathbb{C}),\,\,\omega^2=0\,\,\,\,\,(\omega,\overline{\omega})>0\}
\end{equation}
be the period domain and let $\mathcal{P}\,:\,\,\mathcal{M}_{K3^{[2]}}\,\rightarrow \Omega$ be the period map, where $\mathcal{P}(X,f)=f(\sigma_X)$ and $\sigma_{X}$ is a symplectic 2-form.\\
The period map is surjective and it is a local isomorphism, moreover a global Torelli theorem was proven by Verbitsky in our and in several other cases, see \cite{mar1}, \cite{huy_tor} and \cite{ver}. Whenever this theorem holds it states that two marked \hk manifolds having the same period are birational.\\
Other useful notions are those of Positive, K\"{a}hler and Birational K\"{a}hler cones:
\begin{defn}
Let $X$ be a Hyperk\"{a}hler manifold and let $\omega$ be a K\"{a}hler class.
Let 
\begin{equation}\nonumber
C^+_X=\{l\,\in\,H^{1,1}_{\mathbb{R}}(X),\,\,(l,l)_X>0\}
\end{equation}
be the set of positive classes in $H^{1,1}_{\mathbb{R}}(X)$ and let the positive cone $\mathcal{C}_X$ be its connected component containing $\omega$.\\
Let the K\"{a}hler cone $\mathcal{K}_X\,\subset\,\mathcal{C}_X$ be the set of K\"{a}hler classes.\\
The birational K\"{a}hler cone is the union
\begin{equation}
\mathcal{BK}_X\,=\,\bigcup_{f\,:\,X\,\dashrightarrow\,X'}\,f^*\mathcal{K}_{X'},
\end{equation}
where $f\,:\,X\,\dashrightarrow\,X'$ runs through all birational maps $X\,\dashrightarrow\,X'$ from $X$ to another Hyperk\"{a}hler manifold $X'$.
\end{defn}
There are several results on the structure of these cones, see \cite[Section 27 and 28]{huy} and \cite[Section 9]{mar2}.\\
Last but not least is a result due to Huybrechts on birational maps of Hyperk\"{a}hler varieties (see \cite[Theorem 3.2]{mar1}):
\begin{lem}\label{lem:graph}
Let $(X,f)$ and $(Y,g)$ be two marked Hyperk\"{a}hler manifolds such that $\mathcal{P}(X,f)=\mathcal{P}(Y,g)$ and $X$ and $Y$ are birational. Then there exists an effective cycle $\Gamma=Z+\sum_{j}Y_j$ in $X\times Y$ satisfying the following conditions:
\begin{itemize}
\item Z is the graph of a bimeromorphic map from $X$ to $Y$.
\item The composition $g^{-1}\circ f$ is equal to $\Gamma_*\,:\,H^2(X,\mathbb{Z})\,\rightarrow\,H^2(Y,\mathbb{Z})$.
\item The codimensions of the projections $\pi_1(Y_j)$ and $\pi_2(Y_j)$ are equal.
\item If $\pi_i(Y_j)$ has codimension 1 then it is supported by an effective uniruled divisor.
\end{itemize}
\end{lem}

\subsection{Lattices and discriminant forms}\label{ssec:lattice_theory}
In this subsection we summarize several notions on lattice theory and we analyze some lattices appearing in the rest of the paper. Most of these results are taken from \cite{nik2}.\\
First of all let us start with the basic notions of discriminant groups and forms: given an even lattice $N$ with quadratic form $q$ we can consider the group $A_N=N^{\vee}/N$, which is called discriminant group and whose elements are denoted $[x]$ for $x\in N^\vee$. We denote with $l(A_N)$ the least number of generators of $A_N$. On $A_N$ there is a well defined quadratic form $q_{A_N}$ taking values inside $\mathbb{Q}/2\mathbb{Z}$, which is called discriminant form. Moreover we call $(n_+,n_-)$ the signature of $q$ and therefore of $N$ as a lattice.

\begin{lem}\cite[Corollary 1.13.5]{nik2}\label{lem:nik_spezza}
Let $S$ be an even lattice of signature $(t_+,t_-)$. Then the following hold:
\begin{itemize}
\item If $t_+>0$, $t_->0$ and $t_++t_->2+l(A_S)$ then $S\cong U\oplus T$ for some lattice $T$.
\item If $t_+>0$, $t_->7$ and $t_++t_->8+l(A_S)$ then $S\cong E_8(-1)\oplus T$ for some lattice $T$.
\end{itemize}
\end{lem}

\begin{lem}\cite[Proposition 1.4.1]{nik2}\label{lem:nik_overlattice}
Let $S$ be an even lattice. There exists a bijection $S'\rightarrow H_{S'}$ between even overlattices of finite index of $S$ and isotropic subgroups of $A_S$. Moreover the following hold:
\begin{enumerate}
\item $A_{S'}=(H_{S'}^\perp)/H_{S'}\subset A_S$.
\item $q_{A_{S'}}=q_{A_S|A_{S'}}$.
\end{enumerate}
\end{lem}
We will often need to analyze primitive embeddings of an even lattice into another one. Let us make some useful remarks whose proofs can also be found in \cite{nik2}:
\begin{oss}\label{oss:pre_prim}
A primitive embedding of an even lattice $S$ into an even lattice $N$ is equivalent to giving $N$ as an overlattice of $S\oplus S^{\perp_N}$ corresponding to an isotropic subgroup $H_S$ of $A_S\oplus A_{S^{\perp_N}}$. Moreover there exists an isometry $\gamma\,:\,p_S(H_S)\,\rightarrow\,p_{S^{\perp_N}}(H_S)$ between $q_S$ and $q_{S^{\perp_N}}$ ($p_S$ denotes the natural projection $A_S\oplus A_{S^{\perp_N}}\rightarrow A_S$). Notice that this implies $H_S=\Gamma_\gamma(p_S(H_S))$, where $\Gamma_\gamma$ is the pushout of $\gamma$ in $A_S\oplus A_{S^{\perp_N}}$, \ie
\begin{equation}
\Gamma_\gamma=\{(a,\gamma(a)),\,a\,\in\,H_S\}.
\end{equation}
\end{oss}
\begin{oss}
Suppose we have an even lattice $S$ with signature $(s_+,s_-)$ and discriminant form $q(A_S)$ primitively embedded into an even lattice $N$ with signature $(n_+,n_-)$ and discriminant form $q(A_N)$. Let $K$ be an even lattice, unique in its genus and such that $O(K)\rightarrow O(q_{A_K})$ is surjective, with signature $(k_+,k_-)$ and discriminant form $-q(A_N)$. It follows from \cite{nik2} that primitive embeddings of $S$ into $N$ are equivalent to primitive embeddings of $S\oplus K$ into an unimodular lattice $T$ of signature $(n_++k_+,n_-+k_-)$, such that both $S$ and $K$ are primitively embedded in $T$. By \Ref{oss}{pre_prim} an embedding of $S\oplus K$ into a finite overlattice $V$ such that both $S$ and $K$ are primitively embedded into it is equivalent to giving subgroups $H_S$ of $A_S$ and $H_N$ of $A_N$ and an isometry $\gamma\,:\,q_{A_S|H_S}\,\rightarrow\,-q_{A_N|H_N}$. Finally a primitive embedding of $V$ into $T$ is given by the existence of a lattice with signature $(v_-,v_+)$ and discriminant form $-q_V$. 
\end{oss}
Keeping the same notation as before we give a converse to these remarks:
\begin{lem}\cite[Proposition 1.15.1]{nik2}\label{lem:nik_immerge}
Primitive embeddings of S into an even lattice $N$ are determined by the sets $(H_S,H_N,\gamma,K,\gamma_K)$,where $K$ is an even lattice with signature  $(n_+-s_+,n_--s_-)$ and discriminant form $-\delta$, $\delta\,\cong\,(q_{A_S}\oplus -q_{A_N})_{|\Gamma_\gamma^\perp/\Gamma_\gamma}$ 
and $\gamma_K\,:\,q_K\,\rightarrow\,(-\delta)$ is an isometry.\\
Moreover two such sets $(H_S,H_N,\gamma,K,\gamma_K)$ and $(H'_S,H'_N,\gamma',K',\gamma'_K)$ determine isometric sublattices if and only if
\begin{itemize}
\item $H_S=\lambda H'_S$, $\lambda\in O(q_S)$.
\item There exist $\epsilon\,\in\,O(q_{A_N})$ and $\psi\,\in\,Isom(K,K')$ such that $\gamma'=\epsilon\circ\gamma$ and $\overline{\epsilon}\circ\gamma_K=\gamma'_K\circ\overline{\psi}$. Here $\overline{\epsilon}$ and $\overline{\psi}$ are the isometries induced among discriminant groups.
\end{itemize} 
\end{lem}
The following is a lemma concerning primitive vectors, we include it here since it is needed in the proof of \Ref{lem}{emblemma}:
\begin{lem}\cite[Lemma 7.5]{ghs}\label{lem:ghs_orbit}
Let $T$ be an even lattice such that $T\cong U^2\oplus N$ for some lattice $N$, and let $v,w\in T$ be two primitive vectors such that the following hold:
\begin{itemize}
\item $v^2=w^2$.
\item $(v,T)\cong m\mathbb{Z}\cong (w,T)$, where $(v,T)$ is the image of the linear function $(v,-)$ applied to $T$. 
\item $[\frac{v}{m}]=[\frac{w}{m}]$ in $A_T$.
\end{itemize}
Then there exists an isometry $g$ of $T$ such that $g(v)=w$.
\end{lem}

\begin{ex}\label{ex:discrinvol}
Let us specialize to the lattices of interest to us. Let $E_8(-2)$ be the lattice $E_8$ with quadratic form multiplied by -2. Let $L$ be as in \eqref{latticeL} and let 
\begin{eqnarray}\label{latticeLam}
\Lambda=U^4\oplus E_8(-1)^2.\\\label{latticeM}
M=E_8(-2) \oplus U^3\oplus\,(-2).
\end{eqnarray}
The lattice $E_8(-2)$ has discriminant group $(\mathbb{Z}_{/(2)})^8$ and discriminant form $q_{E_8(-2)}$ given by the following matrix:
\begin{equation}\nonumber
\left(\begin{array}{cccccccc}
1&0&0&\frac{1}{2}&0&0&0&0\\
0&1&\frac{1}{2}&0&0&0&0&0\\
0&\frac{1}{2}&1&\frac{1}{2}&0&0&0&0\\
\frac{1}{2}&0&\frac{1}{2}&1&\frac{1}{2}&0&0&0\\
0&0&0&\frac{1}{2}&1&\frac{1}{2}&0&0\\
0&0&0&0&\frac{1}{2}&1&\frac{1}{2}&0\\
0&0&0&0&0&\frac{1}{2}&1&\frac{1}{2}\\
0&0&0&0&0&0&\frac{1}{2}&1
\end{array}\right).
\end{equation}
The lattice $L$ satisfies
\begin{equation}
A_{L}=\mathbb{Z}_{/(2)},\,\,\,\,\,\,\,\,q_{A_L}(1)=-\frac{1}{2}.
\end{equation}
Since $\Lambda$ is unimodular $A_{\Lambda}=\{0\}$. 
The lattice $(-2)$ has discriminant group $\mathbb{Z}_{/(2)}$ and discriminant form $q'$ with $q'(1)=q_{A_{(-2)}}(1)=\frac{1}{2}$.
Therefore the lattice $M$ has discriminant form $q_{E_8(-2)}\oplus q'$ over the group $(\mathbb{Z}_{/(2)})^9$.
\end{ex}

\begin{lem}\label{lem:biggen_g}
Let $g\in O(L)$, then there exists an embedding $L\subset\Lambda$ and an isometry $\overline{g}\,\in\,O(\Lambda)$ such that $\overline{g}_{|L}=g$ and $\overline{g}_{|L^\perp}=Id$.
\begin{proof}
The isometry $g$ induces an automorphism of the discriminant group $A_L$. Since $A_L=\mathbb{Z}_{/(2)}$ this automorphism is the identity. Let $[v/2]$ be a generator of $A_L$ such that $v^2=-2$. We then have $g([v/2])=[v/2]$, \ie $g(v)=v+2w$. Consider now a lattice of rank 1 generated by an element $x$ of square 2, its discriminant group is still $\mathbb{Z}_{/(2)}$ and is generated by $[x/2]$ with discriminant form given by $q(x/2)=1/2$. Notice that $L\oplus\,\mathbb{Z}x$ has an overlattice isometric to $\Lambda$ which is generated by $L$ and $\frac{x+v}{2}$.\\ We now extend $g$ on $L\oplus x$ by imposing $g(x)=x$ and we thus obtain an extension $\overline{g}$ of $g$ to $\Lambda$ defined as follows:
\begin{eqnarray}\nonumber
\overline{g}(e) &= &g(e),\,\,\,\,\forall\,e\,\in\,L,\\ \nonumber
\overline{g}(x)&= &x,\\ \nonumber
\overline{g}(\frac{x+v}{2})&= &\frac{x+g(v)}{2}.
\end{eqnarray}
\end{proof}
\end{lem}

To conclude this section we analyze the behaviour of (-2) vectors inside $L$ and $M$, since they will play a fundamental role in the proof of \Ref{thm}{main_intro}. We will need the following:
\begin{lem}\label{lem:evenembed}
Let $(-2)$ be the usual lattice and let $e$ be one of its generators. Let $L$ and $M$ be as before.
Then the following hold:
\begin{itemize}
\item Up to isometry there is only one primitive embedding $(-2)\hookrightarrow M$ such that $(e,M)=2\mathbb{Z}$ (\ie $e$ is 2-divisible). Moreover $e\oplus e^\perp=M$.
\item Up to isometry there is only one primitive embedding $(-2)\hookrightarrow L$ such that $(e,L)=2\mathbb{Z}$. Moreover $e\oplus e^\perp=L$.
\end{itemize}

Furthermore all other primitive embeddings into $M$ given by $(H_e,H_M,\gamma,K,\gamma_K)$ satisfy the following:
\begin{equation}\label{halfness_eq}
\exists s\,\in\, A_K,\,\,q_{A_K}(s,s)=\pm\frac{1}{2}.
\end{equation}

\begin{proof}

By \Ref{lem}{nik_immerge} we know that the quintuple $(H_e,H_M,\gamma,K,\gamma_K)$ determines  primitive embeddings of $e$ inside $M$ and the quintuple $(H_e,H_L,\gamma,K,\gamma_K)$ provides those into $L$. A direct computation shows that primitive embeddings of $e$ into $L$ are 2-divisible only for the quintuple $(\mathbb{Z}_{/(2)},A_L,Id,U^3\oplus E_8(-1)^2,Id)$.\\ Now let us move on to the case of $M$. If $H_e=Id$ we have $K\cong U^2\oplus E_8(-2)\oplus (2)\oplus (-2)$, obviously $e$ is not 2-divisible in this case and this satisfies \eqref{halfness_eq}. If $H_e=\mathbb{Z}_{/(2)}$ and $(H_M,A_M^{\perp_{A_{E_8}}})\neq 0$ we obtain nonetheless condition \eqref{halfness_eq}, and again $e$ is not 2-divisible in this embedding since $e\oplus e^{\perp_M}$ is properly contained in $M$ with index a multiple of 2. Therefore $(\mathbb{Z}_{/(2)},A_M^{\perp_{A_{E_8}}},Id,U^3\oplus E_8(-2),Id)$ is the only possible case.

\end{proof}
\end{lem}

\section{Action of automorphisms on cohomology}\label{sec:gen_lattices}

In this section we provide a series of useful facts about finite groups acting faithfully on $X$ and provide a generalization of some results contained in \cite{nik1}. We wish to remark that some among these results are already contained in \cite{beau2}, such as most of \Ref{lem}{gactionlemma} and \eqref{exactgroup}.\\
\begin{defn}\label{defn:inv_locus}
Let $G$ be a finite group acting faithfully on $X$, we define the invariant locus $T_{G}(X)$ inside $H^2(X,\mathbb{Z})$ to be the fixed locus of the induced action of $G$ on cohomology. Moreover we define the co-invariant locus $S_{G}(X)$ as $T_{G}(X)^{\perp}$. 
The fixed locus of $G$ on $X$ will be denoted as $X^{G}$.
\end{defn}
Furthermore if a group $G$ acts on a lattice $R$ we define $T_G(R)$ to be the invariant sublattice and $S_G(R)=T_G(R)^\perp$. 
From now on we keep the same notations as in \cite{nik1}, apart for the following:
\begin{defn}\label{defn:transcendant}
$T(X)$ is the least integer Hodge structure (\ie $T(X)$ is a lattice and $T(X)\otimes\mathbb{C}$ is a Hodge structure) such that $\sigma\,\in\,T(X)\otimes\,\mathbb{C}$ and it is called the transcendant lattice. Furthermore let $S(X)=T(X)^\perp$.
\end{defn}
Let us remark that $S(X)=H^{1,1}_{\mathbb{Z}}(X)$ if $X$ is projective or generic, in fact in those cases $H^{1,1}_{\mathbb{Z}}(X)^\perp$ is an irreducible Hodge structure and $H^{1,1}_{\mathbb{Z}}(X)^\perp\otimes\mathbb{C}$ contains $\sigma$.
This definition differs from that given in \cite{nik1} and in several other papers (where usually $S(X)=H^{1,1}_{\mathbb{Z}}(X)$ and $T(X)=S(X)^\perp$). 
The same definition can be given for any symplectic manifold. 
\begin{ex}
An example where our definition differs from the usual one is given by a very general elliptic K3, where we have $H^{1,1}_{\mathbb{Z}}(X)=\mathbb{Z}v,$ $v^2=0$ and $T(X)\,\cong\,v^\perp/v\,\cong\,U^2\,\oplus\,E_8(-1)^2$, $S(X)\,\cong\,U$. 
\end{ex}
Let $\gamma_X$ be the following useful map:
\begin{equation}\label{sigmamap}
\gamma_X\,:\,T(X)\,\rightarrow\,\mathbb{C}.
\end{equation}
Given by $\gamma_X(x)=(\sigma,x)_X$, which has kernel $T(X)\cap S(X)=0$.\\
 Moreover we have the following exact sequence for any finite group $G$ of Hodge isometries on $H^2(X,\mathbb{Z})$:
\begin{equation}\label{exactgroup}
1\,\rightarrow\, G_0\,\rightarrow\,G\,\stackrel{\pi}{\rightarrow}\,\Gamma_m\,\rightarrow\,1,
\end{equation}
where $\Gamma_m\subset U(1)$ is a cyclic group of order $m$. In fact the action of $G$ on $H^{2,0}$ is the action of a finite group on $\mathbb{C}^*$. 
We also denoted $G_0=ker(G\,\rightarrow\,Aut(H^{2,0}(X)))$.
  The following result is a generalization of \cite[theorem 3.1]{nik1}:
\begin{lem}\label{lem:gactionlemma}
Let $X$ be a Hyperk\"{a}hler 4-fold of $K3^{[2]}$-type and $G\,\subset\,Aut(X)$, $|G|$ finite. Then the following hold:
\begin{enumerate}
\item $g\in G$ acts trivially on $T(X)\,\iff\,g\in G_0$.
\item The representation of $\Gamma_m$ on $T(X)\otimes\mathbb{Q}$ splits as the direct sum of irreducible representations of the cyclic group $\Gamma_m$ having maximal rank (\ie of rank $\phi(m)$).
\end{enumerate}
\begin{proof}
\begin{enumerate}
\item Let $g\in G_0$. Let us show that $g^*$ acts trivially on $T(X)\otimes\mathbb{Q}$. We start by considering the kernel of the map $g^*-Id_{T(X)}$ which is a lattice (and a Hodge substructure) $R$ inside $T(X)$. Hence, by minimality of $T(X)$, $R\otimes\mathbb{Q}$ is either 0 or $R\otimes\mathbb{Q}=T(X)\otimes\mathbb{Q}$. Considering the map \eqref{sigmamap}, since $g^*$ is a Hodge isometry we have
\begin{equation}\nonumber
\gamma_X(x)=(g^*\sigma,g^*x)=(\sigma,g^*x).
\end{equation}
Since $g^*\sigma=\sigma$ we have that $g^*x-x\,\in\,ker(\gamma_X)=T(X)\cap S(X)=0$. Thus $R$ is all of $T(X)$.\\ 
To obtain the converse we prove that $g^*\sigma=\lambda\sigma$ with $\lambda\neq1$ implies that 1 is not an eigenvalue of $g^*$ on $T(X)$. In fact
\begin{equation}
\gamma_X(x)=(g^*\sigma,g^*x)=\lambda\gamma_X(g^*x), \nonumber
\end{equation}
\ie $g^*x\neq x$.
\item The preceeding arguments show that every nontrivial element of $G/G_0$ has no eigenvalue 1 on $T(X)$ and hence also on $T(X)\otimes\mathbb{Q}$, this implies our claim.
\end{enumerate}
\end{proof}
\end{lem}

Let now $G$ be a finite group of automorphisms such that $G=G_0$. Following Nikulin we will call such $G$ an algebraic automorphism group. We want to give some useful generalizations of \cite[section 4]{nik1}:
\begin{lem}\label{lem:algaction}
Let G be a finite algebraic automorphism group of a fourfold $X$ of $K3^{[2]}$-type, then the following assertions are true:
\begin{enumerate}
\item $S_G(X)$ is nondegenerate and negative definite.
\item $S_G(X)$ contains no element with square -2.
\item $T(X)\subset T_G(X)$ and $S_G(X)\subset S(X)$.
\item $G$ acts trivially on $A_{S_G(X)}$.
\end{enumerate}
\begin{proof}
The third assertion is an immediate consequence of \Ref{lem}{gactionlemma} because $G$ acts as the identity on $\sigma$ and therefore on all of $T(X)$.\\
To prove that $S_G(X)$ and $T_G(X)$ are nondegenerate let $H^2(X,\mathbb{Z})=\oplus_{\rho}U_\rho$ be the decomposition in orthogonal representations of $G$, where $U_{\rho}$ contains all irreducible representations of $G$ of character $\rho$ inside $H^2(X,\mathbb{Z})$. Obviously $T_G(X)=U_{Id}$ and $S_G(X)=\oplus_{\rho\neq Id}U_\rho$, which implies they are orthogonal and of trivial intersection. Hence they are both nondegenerate.\\
 Since $G$ is finite there exists a $G$-invariant K\"{a}hler class $\omega_G$ given by $\sum_{g\in G}g\omega$, where $\omega$ is any K\"{a}hler class on $X$.  
Therefore we have: 
\begin{equation}\nonumber
\sigma\mathbb{C}\oplus\overline{\sigma}\mathbb{C}\oplus\omega_G\mathbb{C}\,\subset\,T_G(X)\otimes\mathbb{C}.
\end{equation}
Hence the lattice $S_G(X)$ is negative definite.\\
To prove the last assertion let us proceed as in \Ref{lem}{biggen_g}, \ie let us choose a primitive embedding of $H^2(X,\mathbb{Z})$ in the lattice $\Lambda$ such that the action of $G$ extends trivially outside the image of $H^2(X,\mathbb{Z})$. Therefore $S_G(X)\cong S_G(\Lambda)$ and $A_{S_G(\Lambda)}\cong A_{T_G(\Lambda)}$, where the isomorphism is $G$ equivariant. $G$ acts trivially on $T_G(\Lambda)$, thus its induced action on $A_{T_G(\Lambda)}$ is trivial. Using the $G$ equivariant isomorphism we have that $G$ acts trivially also on $A_{S_G(\Lambda)}=A_{S_G(X)}$.\\
Let us prove that there are no $-2$ vectors inside $S_G(X)$. Assume on the contrary that we have an element $c\in S_G(X)$ such that $(c,c)=-2$. Then by \cite[Theorem 1.12]{mar2} it is known that either $\pm c$ or $\pm 2c$ is represented by an effective divisor D on X. Let $D'=\sum_{g\in G}gD$ which is also an effective divisor on X, but $[D']\in\,S_G(X)\cap T_G(X)=\{0\}$. This implies $D'$ is linearly equivalent to 0, which is impossible.
\end{proof}
\end{lem}

Now we can use \Ref{lem}{graph} to give sufficient conditions for an isometry  $\psi$ of $L$ to be induced by a birational map $\psi'$ of some marked Hyperk\"{a}hler manifold $(X,f)$ such that $f\circ\psi'^*\circ f^{-1}=\psi$. Thus we obtain a generalization of \cite[Theorem 4.3]{nik1}: 
\begin{thm}\label{thm:cohom_to_aut}
Let $G$ be a finite subgroup of $O(L)$. Suppose that the following hold:
\begin{enumerate}
\item $S_G(L)$ is nondegenerate and negative definite.
\item $S_G(L)$ contains no element with square $(-2)$.
\end{enumerate}
Then $G$ is induced by a subgroup of $Bir(X)$ for some manifold $(X,f)$ of $K3^{[2]}$-type.
\begin{proof}
By the surjectivity of the period map and by \Ref{lem}{algaction} we can consider a marked $K3^{[2]}$-type 4-fold $(X,f)$ such that $T(X)\stackrel{f}{\rightarrow} T_G(L)$ is an isomorphism and also $S(X)\stackrel{f}{\rightarrow} S_G(L)$ is.\\
Let $g\in G$, let us consider the marked varieties $(X,f)$ and $(X,g\circ f)$. They have the same period in $\Omega$ and hence by \Ref{lem}{graph} we have $f^{-1}\circ g\circ f=\Gamma_*$. Here $\Gamma=Z+\sum_{j}Y_j$ in $X\times X$, where $Z$ is the graph of a bimeromorphic map from $X$ to itself and $Y_j$'s are cycles with $codim(\pi_i(Y_j))\geq 1$.\\ 
 We will prove that all $Y_j$'s contained in $\Gamma$ have $codim(\pi_i(Y_j))>1$, thus implying $\Gamma_*=Z_*$ on $H^2_{\mathbb{Z}}$. We know those of codimension 1 are uniruled and effective, moreover it is known (see \cite[Proposition 28.7]{huy}) that uniruled divisors cut out the closure of the birational K\"{a}hler cone $\mathcal{BK}_X$, \ie $(\alpha,D)\geq0$ for all $\alpha\in\overline{\mathcal{BK}}_X$ and for all uniruled $D$. We wish to remark that the manifold $X$ we chose has $\overline{\mathcal{BK}}_X=\overline{\mathcal{C}}_X$ by \cite[Theorem 9.17]{mar1} (it contains no -2 divisors).\\ Let $\beta\in\mathcal{C}_X$ be a K\"{a}hler class and let $D\in Pic(X)$ be a uniruled divisor, we can write
\begin{equation}\nonumber
\beta=\alpha+\gamma,\,\,f(\alpha)\,\in\,T_G(L)\otimes\mathbb{R},\,\,f(\gamma)\,\in\,S_G(L)\otimes\mathbb{R}.
\end{equation}
Hence $0<(\beta,D)=(\gamma,D)$ and moreover we have $(f^{-1}\circ g\circ f(\beta),D)=(f^{-1}\circ g\circ f(\gamma),D)=(\gamma,f^{-1}\circ g^{-1}\circ f(D))\geq0$ because $f^{-1}\circ g\circ f(\beta)\in \overline{\mathcal{BK}}_X$ and $D$ is uniruled. Here is the contradiction: 
\begin{equation}\nonumber
0<(\beta,\sum_{h\in G}f^{-1}\circ h\circ f(D)),
\end{equation}
which implies $0\neq D'=\sum_{h\in G}hD\in f^{-1}(T_G(L)\cap S_G(L))=0$, hence there are no uniruled divisors inside $Pic(X)$. Moreover we obtain $\Gamma_*=Z_*$, \ie there exists a bimeromorphic map $\psi'$ of $X$ such that $\psi'^*=f^{-1}\circ g\circ f$ on $H^2(X)$. 
\end{proof}
\end{thm}

\section{Fixed locus of a Symplectic involution}
Our work on symplectic involutions starts with an analysis of  the fixed locus of a symplectic involution $\varphi$ on $X$. The main result of this section is the following: 
\begin{thm}\label{thm:fixedpoints}
Let $X$ be a Hyperk\"{a}hler manifold of $K3^{[2]}$-type with a symplectic involution $\varphi$. Then \Ref{conj}{camere_c} holds true, the fixed locus $X^{\varphi}$ consists of 28 isolated points and one K3 surface. 
Moreover the lattice $T_{\varphi}(X)$ has rank 15.
\end{thm}
Notice that this is what happens in \Ref{ex}{s2}.
The starting point of our proof will be the following result of Camere:
\begin{prop}\cite{cam}\label{prop:camerefixedpoints}
Let $X$ be a manifold of \ktipo and let $\varphi\subset\,Aut(X)$ be a symplectic involution. Then $rank(T_{\varphi}(X))\geq 11$. Moreover, unless $X^{\varphi}$ contains a complex torus, we have $rank(T_{\varphi}(X))=15$ and $X^{\varphi}$ consists of 28 isolated points and a K3 surface.  
\end{prop}

\begin{oss}\label{oss:noidentity}
\Ref{prop}{camerefixedpoints} implies that a symplectic involution on $X$ cannot induce the identity map on cohomology. The same result can be proven for automorphisms of manifolds of $K3^{[2]}$-type of any order.
\end{oss}

\begin{ex}\label{ex:s2}
Let $S$ be a K3 surface endowed with a symplectic involution $\psi$. Then the manifold $X=S^{[2]}$ has a symplectic involution $\psi^{[2]}$ fixing 28 points and 1 K3 surface $Y$. Notice that $Y$ is the minimal resolution of $S_/\psi$.\\
For further examples the reader can consult \cite{cam}.
\end{ex}
To prove \Ref{thm}{fixedpoints} let us do some preliminary work to analyze small deformations of  the couple $(X,\varphi)$, \ie the following:
let us choose a small ball $U$ representing $Def(X)$, whose tangent space at the origin is given by $H^1(\mathcal{T}_X)$. The symplectic involution on $X$ extends to an automorphism of the versal deformation family $\mathcal{X}\,\rightarrow\, U$ as follows:
\begin{equation}\nonumber
\begin{array}{ccc}
\mathcal{X} & \stackrel{M}{\longrightarrow} & \mathcal{X}\\
\downarrow &  & \downarrow\\
U & \stackrel{m}{\longrightarrow} & U.
\end{array}
\end{equation}
Here $m$ is the involution on $U$ (which is small enough to have $m(U)=U$), induced by the action of the symplectic involution $\varphi$ on $H^1(\mathcal{T}_X)$. Moreover $m$ induces an involution of $\mathcal{X}$ which yields fibrewise isomorphisms between $\mathcal{X}_t$ and $\mathcal{X}_{m(t)}$. The differential of $m$ at $0$ is given by the action of $\varphi$ on $H^1(\mathcal{T}_X)$, which is the same as the action on $H^{1,1}(X)$ since the symplectic form $\sigma$ induces an isomorphism between those two and $\sigma$ is preserved by the action of $\varphi$. On the other hand $U^m$ is smooth, since $m$ is linearizable, and hence 
\begin{equation}\nonumber
dim (U^m)=rank (T_\varphi(X))-2,
\end{equation}
which is always positive by \Ref{lem}{algaction}. We wish to obtain a deformation of the couple $(X,\varphi)$, hence we need to restrict to $U^m$ to get a fibrewise involution. Therefore we obtain the following diagram:
\begin{equation}\label{deform}
\begin{array}{cccc}
\mathcal{Y}= &\mathcal{X}_{|U^m} & \stackrel{M}{\longrightarrow} & \mathcal{X}\\
&\downarrow &  & \downarrow\\
&U^m & \stackrel{m}{\longrightarrow} & U,
\end{array}
\end{equation}
where $\mathcal{Y}\,\rightarrow\, U^m$ represents the functor of deformations of the couple $(X,\varphi)$, \ie  all small deformations of this couple must embed in $\mathcal{Y}\,\rightarrow\, U^m$. The involutions $\varphi_{t}$ are given by $M_{|X_t}$.\\
It is obvious that this deformation space is "maximal" in some sense. Let us make this more precise using the period map.

\begin{defn}
Given a finite group $G\subset O(L)$ we denote $\Omega_G$ the set of points $(X,f)$ in the period domain such that $f(\sigma_X)\,\in\,T_G(L)$. 
\end{defn}
\begin{defn}
Given $(X,f)$ with a group $G$ acting faithfully on it via symplectic bimeromorphic maps, we call the following a maximal family of deformations of $(X,G_{Bir})$
\begin{equation}\nonumber
\begin{array}{ccc}
X & \stackrel{i}{\longrightarrow} & \mathcal{X}_U\\
\downarrow &  & \downarrow\\
\{0\} & \stackrel{i}{\longrightarrow} & U,
\end{array}
\end{equation}
where the family $\mathcal{X}$ over $U$ is endowed with a fibrewise faithful bimeromorphic action of $G$ and the period map $\mathcal{P}$, given a compatible marking, sends surjectively a neighbourhood of $0\in\,U$ inside a neighbourhood of $\mathcal{P}(X,f)\cap \Omega_{G}$.\\
We give the same definition for maximal families $(X,G_{Aut})$ or $(X,G_{Hod})$ having $G$ acting as symplectic automorphisms or Hodge isometries on $H^2(X,\mathbb{Z})$ respectively. 
Notice that the family $\mathcal{Y}\,\rightarrow\,U^m$ we stated before is a maximal family for the couple $(X,\varphi)$.

\end{defn}
\begin{oss}\label{oss:generic_point}
We remark that the set $\Omega'_G=\bigcup_{v\in T_G(L)}\{x\in \Omega_G\,:\,(x,v)=0\}$ is the union of countable codimension 1 subsets and consists of Hodge structures on marked varieties $(X,f)$ over $\Omega_G$ such that the inclusion\\ $f(T(X))\hookrightarrow T_G(L)$ is proper. Moreover outside this set $T(X)$ is irreducible. 
\end{oss}

Now we can use this construction to prove the following fact: 
\begin{prop}\label{prop:fixed_abel}
Let $(X,\varphi)$ be as before and suppose $\varphi$ fixes at least one complex torus $T$. Then $T_{\varphi}(X)$ has rank at most 6.
\begin{proof}
Suppose on the contrary that $T_{\varphi}(X)$ has rank $\geq7$. Let us consider small deformations of the couple $(X,\varphi)$ over a representative $U$ of $Def(X)$ given by
\begin{equation}
\begin{array}{ccc}
\mathcal{X}_{|U^m} & \stackrel{\Phi}{\longrightarrow} & \mathcal{X}\\
\downarrow &  & \downarrow\\
U^m & \stackrel{m}{\longrightarrow} & U.
\end{array}
\end{equation}
As shown in \eqref{deform}. We let $\sigma_t$ be the symplectic form on $\mathcal{X}_t$.\\ We remark that, by linear algebra, the fixed locus $X^{\varphi}$ is smooth and consists only of symplectic varieties since the symplectic form $\sigma$ restricts to a nonzero symplectic form on all connected components of $X^\varphi$. Moreover it is stable for small deformations of the couple $(X,\varphi)$, \ie the fixed locus $\mathcal{X}^{\Phi}$ is a small deformation of the fixed locus $X^{\varphi}$.
Therefore we have a well defined map of integral Hodge structures $H^2(\mathcal{X}_t,\mathbb{C})^{\Phi_t}\,\rightarrow\, H^2(T_t,\mathbb{C})$ sending a class on $H^2(\mathcal{X}_t)$ to its restriction to $T_t$, where $T_t$ is a small deformation of $T$ fixed by $\Phi_t$ (\ie is a component of the fibre over $t$ of $\mathcal{X}^{\Phi}$). Since $\Phi_t(\sigma_t)=\sigma_t$ and $\sigma_{t|T_t}\neq 0$ this map is not the zero map and, being a map of Hodge structures, its kernel is again a Hodge structure.\\ 
Given a marking $F$ over $\mathcal{X}$ we have that $(\mathcal{X},F)$ is a maximal family of deformations of the couple $(X,\varphi)$. Let $V=\{\mathcal{P}(\mathcal{X}_t,F_t),\,t\in U\}\,\subset\,\Omega_\varphi$, by \Ref{oss}{generic_point} there exists $u\in V\backslash \Omega'_\varphi$ and this period corresponds to a marked manifold $(\mathcal{X}_t,F_t)$ such that $T(\mathcal{X}_{t})=T_{\Phi_t}(\mathcal{X}_{t})$, \ie this Hodge structure is irreducible. Therefore we have that the map $H^2(\mathcal{X}_t,\mathbb{C})^{\Phi_t}\rightarrow H^2(T_t,\mathbb{C})$ is an injection. But this is absurd if $T_\varphi(X)$ has rank greater than 6 since $H^2(T_t)$ has dimension 6.
\end{proof}
\end{prop}
\begin{dim2}\nonumber
By \Ref{prop}{camerefixedpoints} we have that $rank(T_{\varphi}(X))\geq 11$. By \Ref{prop}{fixed_abel} we therefore have that symplectic involutions cannot fix complex tori, hence we have our claim.
\end{dim2}

\section{Deformation equivalence of couples $(X,\varphi)$}

Having determined $X^\varphi$ we proceed to compute $S_{\varphi}(X)$. We will use part of Nikulin's theorems summarized in \Ref{ssec}{lattice_theory}, and also the following result concerning invariant and co-invariant lattices of involutions on $L$ and $\Lambda$:
\begin{lem}\label{lem:Lam_lem}
\begin{enumerate}
Let $\varphi\in O(L)$ be an involution and let $L\,\subset \Lambda$ as in \Ref{lem}{biggen_g}, \ie $\varphi\subset O(L)\subset O(\Lambda)$. Then the following hold:
\item The quotient $\Lambda/(T_{\varphi}(\Lambda)\,\oplus\, S_{\varphi}(\Lambda))$ is of $2-$torsion. 
\item $S_{\varphi}(L) \cong S_{\varphi}(\Lambda)$.
\end{enumerate}
\begin{proof}
Given an element $t\,\in\,\Lambda$ we have $t=\frac{t+\varphi(t)}{2}+\frac{t-\varphi(t)}{2}$ and clearly $t+\varphi(t)\,\in\,T_{\varphi}(\Lambda)$ while $t-\varphi(t)\,\in\,S_{\varphi}(\Lambda)$ so we have $2t\,\in\,S_{\varphi}(\Lambda)\,\oplus\,T_{\varphi}(\Lambda)$. The lattice $\Lambda$ is generated by $L$ and $\frac{x+v}{2}$ as in \Ref{lem}{biggen_g}, where $<x>=L^\perp$ and $v$ is a 2-divisible vector of square $-2$ inside $L$.\\ Moreover a vector $b=a\frac{x+v}{2}+w'$, $w'\in L$ is in $S_{\varphi}(\Lambda)$ if and only if
\begin{equation}
-a\frac{x+v}{2}-w'=\varphi(a\frac{x+v}{2}+w')=a\frac{x}{2}+\varphi(w')+a\frac{v}{2}+aw. \nonumber
\end{equation}
Here $\varphi(v)=v+2w$.
But this can happen only if $a=0$, \ie $b\in S_{\varphi}(L)$ 
\end{proof}
\end{lem}
\begin{thm}\label{thm:invol_lattice}
Let $X$ be a Hyperk\"{a}hler manifold of $K3^{[2]}$-type and $\varphi\in Aut(X)$ a symplectic involution. Then the lattice $S_{\varphi}(X)$ is isomorphic to $E_8(-2)$ and $T_{\varphi}(X)$ is isomorphic to $E_8(-2)\,\oplus\, U^3\,\oplus\, (-2)$.

\begin{proof}
We know that $S_{\varphi}(X)$ has rank 8 by \Ref{thm}{fixedpoints} and it equals $S_{\varphi}(\Lambda)$ by \Ref{lem}{Lam_lem}, therefore the discriminant group $A_{S_{\varphi}(\Lambda)}$ can be generated by 8 elements and so does its unimodular complement $A_{T_{\varphi}(\Lambda)}$. This means that we can apply \Ref{lem}{nik_spezza} obtaining $T_{\varphi}(\Lambda)=U\oplus T'$, which means that we can define an involution of $U^3\oplus E_8(-1)^2$ having $S_{\varphi}(\Lambda)$ as the anti-invariant lattice. By \Ref{lem}{algaction} this involution satisfies the conditions of \cite[Theorem 4.3]{nik1} which implies that this involution on $U^3\oplus E_8(-1)^2$ is induced by a symplectic involution $\psi$ on some $K3$ surface $S$ and hence also $S_{\psi}(S)\,\cong\,S_{\varphi}(X)$.\\ Thus, by the work of Morrison on involutions \cite{mor}, we know $S_{\varphi}(X)=E_8(-2)$ and $T_{\varphi}(X)$ is just its orthogonal complement in $L$, which is easily proven to be $E_8(-2)\,\oplus\, U^3\,\oplus\, (-2)$ using \Ref{lem}{nik_immerge}.
\end{proof}
\end{thm}

\begin{cor}\label{cor:conjlattice}
Let $M_1,M_2\,\subset L$ such that $M_1\cong M_2\cong E_8(-2)$. Then there exists $f\,\in\,O(L)$ such that $f(M_1)=M_2$. 
\begin{proof}
By \Ref{ex}{discrinvol} we know the discriminant form and group of $E_8(-2)$. Therefore we can apply \Ref{lem}{nik_immerge}, obtaining that embeddings of $E_8(-2)$ into $L$ are given by quintuples $(H,H',\gamma,K,\gamma_K)$. Moreover two such embeddings $(H,H',\gamma,K,\gamma_K)$ and $(N,N',\gamma',K',\gamma'_{K'})$ are conjugate if and only if we have $H$ conjugate to $N$ through an automorphism of $(\mathbb{Z}_{/(2)})^8$ sending $\gamma$ into $\gamma'$. In our case the computations are particularly simple: due to the values of $q_{E_8(-2)}$ (all elements have square $0$ or $1$) and $q'$ (all nonzero elements have square $\frac{1}{2}$) the only possible choices of $H$ and $H'$ are given by the one element group and so we obtain our claim.\\
Moreover this implies that we can always choose a marking of $(X,\varphi)$ such that the induced action of $\varphi$ on $L$ is given by leaving $(-2)\oplus U^3$ invariant and exchanging the two remaining $E_8(-1)$, so that $S_{\varphi}$ is given by the differences $a-\varphi(a)$ for $a\in\,E_8(-1)$.    
\end{proof}
\end{cor}
Now we can proceed to the proof of \Ref{thm}{main_intro}, \ie that all couples $(X,\varphi)$ where $\varphi$ is a symplectic involution are standard. 
Let us start by defining a space containing any $(X,\varphi)$:
\begin{defn}
Let $\mathcal{M}_{2}$ be the subset of $\mathcal{M}_{K3^{[2]}}$ given by the marked manifolds  $(X,f)$ such that:
\begin{equation}\nonumber
\exists\,V\cong E_8(-2),\,V\subset L\,:\,V\,\subset\,f(H^{1,1}_{\mathbb{Z}}(X)).
\end{equation}
\end{defn}
\begin{prop}
Let $(X,\varphi)$ be a couple consisting of a $K3^{[2]}$-type manifold and a symplectic involution of $X$. Then $\mathcal{M}_{2}$ contains all couples $(X,f)$ for any marking $f$. Moreover the generic point of $\mathcal{M}_2$ corresponds to a marked manifold $(Y,g)$ having a bimeromorphic involution.
\begin{proof}
$\mathcal{M}_{2}$ is locally given by 8 linearly independent conditions on the image of the period map, \ie $\mathcal{P}(X,f)\,\perp\,a$ with $a$ ranging through a set of generators for a lattice of type $E_8(-2)\,\subset\,L$. Due to \Ref{cor}{conjlattice} we assume that the marking is fixed. Given such an $(X,f)$ we can define an involution $\varphi$ inside its cohomology by imposing $f(S_{\varphi}(X))=E_8(-2)\subset\,L$.\\ Since this is a maximal family of Hodge involutions, the generic element $(Y,g)$ of this space has $Pic(Y)=E_8(-2)$ by \Ref{oss}{generic_point} and we know by \Ref{thm}{cohom_to_aut} that $\varphi$ extends to a birational involution on $Y$. Finally a couple $(X,\varphi)$ endowed with a marking $f$ satisfies the condition $E_8(-2)\,\subset\,f(Pic(X))$ by \Ref{thm}{invol_lattice} and is thus inside this space.
\end{proof}
\end{prop}
\begin{defn}
Let $\Omega_2=\mathcal{P}(\mathcal{M}_2)$ and furthermore let $\Omega_{v,2}$ denote the set of $\omega\,\in\Omega_2$ such that $(v,\omega)=0$. 
\end{defn}
Let $M$ be as in \eqref{latticeM}, there is a sublattice $M_0$ of $L$ isomorphic to $M$ given by $f(T_\varphi(S^{[2]}))$, where $(S^{[2]},f)$ is a marked Hyperk\"{a}hler manifold and $\varphi$ is a symplectic involution on it. Moreover, by \Ref{cor}{conjlattice}, all such lattices are conjugate through an isometry of $L$, hence without loss of generality we fix $M_0\subset L$, $M\cong M_0$ and we can impose
\begin{equation}\nonumber
\mathcal{P}(X,f)\in\mathbb{P}(M_0\otimes\mathbb{C})
\end{equation}
for all couples $(X,\varphi)$ and an appropriate marking $f$.

\begin{lem}\label{lem:emblemma}
Let $0\neq w\in M$ be a primitive isotropic vector, then there exist a sublattice $w\in T\subset M$ and a $(-2)$ vector $p$ such that:
\begin{itemize}
\item $p$ is 2-divisible in $M$,
\item $q_{M|T}$ is nondegenerate,
\item $R:=T^{\perp_M}\cong U\oplus <p> \oplus R'$ for some lattice $R'$.
\end{itemize} 

\begin{proof}
Since $M=U^2\oplus (U\oplus E_8(-2)\oplus (-2))$ we can apply \Ref{lem}{ghs_orbit}. Therefore we can analyze up to isometry all isotropic vectors inside $M$ knowing only their divisibility $m$ (\ie $(w,M)=m\mathbb{Z}$) and their image $[\frac{w}{m}]$ in $A_M$. Let us give a basis of $M$ as follows:
\begin{equation}
\{e_1,f_1,e_2,f_2,e_3,f_3,a_1,a_2,a_3,a_4,a_5,a_6,a_7,a_8,t\},
\end{equation}
where $\{e_i,f_i\}$ is a standard basis of $U$, $\{a_1,\dots,a_8\}$ is a standard basis of $E_8(-2)$ and $t$ is a generator of the lattice $(-2)$.\\
The first key remark is that since $A_M$ is of 2-torsion $m$ can either be $1$ or $2$.
Therefore if $m=1$ we have that $\frac{w}{m}$ lies in $M$, which implies $[\frac{w}{m}]=0$ in $A_M$. Thus by \Ref{lem}{ghs_orbit} there exists an isometry $g$ of $M$ sending $w$ to $e_1$. To obtain our claim we let $T=g^{-1}(<e_1,f_1>)$, $p=g^{-1}(t)$ and $R=g^{-1}(<e_2,f_2,e_3,f_3,a_1,a_2,a_3,a_4,a_5,a_6,a_7,a_8,t>)$.\\
If $m=2$ we have that $\frac{w}{2}$ is a square zero element of $M^\vee$, \ie $[\frac{w}{2}]$ has square zero in $A_M$. Looking at \Ref{ex}{discrinvol} it is easy to see that square zero elements must lie in $A_{E_8(-2)}\subset A_M$ and they are given by $[\frac{v}{2}]$ where $v$ is a primitive vector of square $c\equiv 0\,\,\,mod\,\,8$ inside $E_8(-2)$. 
Therefore by \Ref{lem}{ghs_orbit} there exists an isometry $g$ of $M$ sending $w$ to $r=2e_1+\frac{c}{4}f_1+v$. Thus we set $T=g^{-1}(<r,f_1>)$, $p=g^{-1}(t)$, $K=v^{\perp_{E_8(-2)}}$ and $R=g^{-1}(<e_2,f_2,e_3,f_3,K,t>)$.   
\end{proof}

\end{lem}
\begin{lem}\label{lem:denselemma}
Let $0\neq w_0\in M_0$ be a primitive vector of square 0 and let $k\in\mathbb{Z}$.\\ There exists a sequence $\{w_n\}$ of non-zero primitive vectors of norm $2k$  such that $[w_n]$ converges to $[w_0]$ in $\mathbb{P}(M_0\otimes \mathbb{R})$. Moreover if $k$ is odd we may assume that for all $n$:
\begin{equation}\nonumber
(w_n,L)=2\mathbb{Z}.
\end{equation}
Furthermore if $k=0$ we may assume that there exists an element $q$ of square $-2$ and divisibility $2$ such that $w_n\perp q$.
\begin{proof}
We keep the same notation as in \Ref{lem}{emblemma} and we fix an isometry $\eta\,: M_0\,\rightarrow\,M$.\\
Let $w=\eta(w_0)$, since it satisfies the hypothesis of \Ref{lem}{emblemma} we have a lattice $p\oplus U$ orthogonal to $w$, where $p$ is a 2-divisible $(-2)$ vector. Let $e$ and $f$ be two standard generators of such $U$. The sequence $\{[\eta^{-1}(nw+e+kf)]\}_n$ converges to $[w_0]$ in $\mathbb{P}(M_0\otimes\mathbb{R})$ and consists of primitive vectors of square $2k$.\\ If $k$ is odd the sequence $\{[\eta^{-1}(2nw+p+2e+(k+1)f)]\}_n$ converges to $[w_0]$ in $\mathbb{P}(M_0\otimes\mathbb{R})$, we need only to prove that it is composed by 2-divisible vectors. Obviously this is equivalent to proving that $\eta^{-1}(p)$ is 2-divisible in $L$, \ie to proving that $\eta^{-1}(p)\oplus \eta^{-1}(p)^{\perp_L}=L$. We know that $p^{\perp_M}\cong U^3\oplus E_8(-2)$ hence $\eta^{-1}(p)^{\perp_L}$ is an overlattice of $U^3\oplus E_8(-2)^2$ which, by \Ref{lem}{evenembed}, implies $\eta^{-1}(p)$ is 2-divisible in $L$.    
To obtain our last claim we let $q=\eta^{-1}(p)$, as we just proved it is $2-$divisible, of square $-2$ and it is orthogonal to $w_n$.
\end{proof}  
\end{lem}
\begin{defn}
Let $\mathcal{P}_{exc}=\{f\,\in\,M_0\,:\,f^2=-2\,\,,\,(f,L)=2\mathbb{Z}\}$ be the set of exceptional primitive classes inside $M_0$.
\end{defn}
Notice that $\mathcal{P}^{-1}(v)$ contains the Hilbert square of a $K3$ surface for all $v\in\,\mathcal{P}_{exc}$.
\begin{lem}\label{lem:denselemma2}
$\cup_{v\in\mathcal{P}_{exc}}\Omega_{v,2}$ is dense in $\Omega_2$.
\begin{proof}
It is enough to prove that $\cup_{v\in\mathcal{P}_{exc}}\Omega_{v,2}$ is dense in $\Omega_2\cap\mathbb{P}(M_0\otimes\mathbb{C})$ by \Ref{cor}{conjlattice}.\\ Let $Q_{M_0}$ be the subset of isotropic vectors inside $\mathbb{P}(M_0\otimes\mathbb{C})$. Let $Q_{M_0}(\mathbb{R})$ and $Q_{M_0}(\mathbb{Q})$ be the subsets of isotropic vectors spanned by real (respectively rational) isotropic vectors. Let $\omega$ be in $\Omega_2\cap\mathbb{P}(M_0\otimes\mathbb{C})$, we have $\omega^{\perp_{M_0}}\,\cap\,Q_{M_0}(\mathbb{R})=(\alpha\omega+\overline{\omega\sigma})^{\perp_{M_0}}\,\cap\,Q_{M_0}(\mathbb{R})$.\\ But since $(\alpha\omega+\overline{\alpha\omega})^{\perp_{M_0}}$ has signature (1,12) we have that $\exists\,\,u\,\in\,Q_{M_0}(\mathbb{R})\cap (\alpha\omega+\overline{\alpha\omega})^{\perp_{M_0}}$. Since $Q_{M_0}(\mathbb{Q})$ is non-empty it is dense inside $Q_{M_0}(\mathbb{R})$, therefore $\exists\,\{v_n\}$ such that $[v_n]\rightarrow [u]$ in $\mathbb{P}(M_0\otimes\mathbb{C})$, where the $v_n$ are primitive isotropic vectors inside $M_0$. Thus we can apply \Ref{lem}{denselemma} to find a sequence $\{w_n\}$ of elements of $\mathcal{P}_{exc}$ and a sequence $\{v'_n\}$ of isotropic primitive vectors such that $[v'_n]\rightarrow [u]$ and $w_n\perp v'_n$.
\end{proof}
\end{lem}
\begin{dim1}\nonumber
Let $f$ be a marking of $X$ such that $\mathcal{P}(X,f)\subset\mathbb{P}(M_0\otimes\mathbb{C})$ and $f(S_\varphi(X))\perp M_0$.  
Moreover let $\mathcal{X}\rightarrow U$ be a maximal family of deformations of the couple $(X,\varphi)$ as in \eqref{deform} and let $F$ be a marking of $\mathcal{X}$ compatible with $f$ such that $V=\{\mathcal{P}(\mathcal{X}_t,F_t),\,t\in U\}$ is a small neighbourhood of $\mathcal{P}(X,f)$. By \Ref{lem}{denselemma2} there exist a point $v\in V$ and a 2-divisible primitive vector $e$ of square $(-2)$ such that $v\perp e$. Since the global Torelli theorem holds we can use \Ref{lem}{graph} on the manifold $\mathcal{X}_u$ such that $\mathcal{P}(\mathcal{X}_u,F_u)=v$. This gives that $\mathcal{X}_u$ is bimeromorphic to the Hilbert square of a certain K3 surface $S$.\\ Thus we get a bimeromorphic involution $\varphi$ on $S^{[2]}$ such that $S_{\varphi}(S^{[2]})\subset Pic(S)\subset Pic(S^{[2]})$, where 
\begin{equation}\nonumber
Pic(S)=\{t\in Pic(S^{[2]}), e\perp t\}.
\end{equation}
By \cite[Theorems 4.3 and 4.7]{nik1} we have a symplectic involution $\psi$ on $S$ given by the action of $\varphi$ on $e^\perp\cong H^2(S,\mathbb{Z})$ which induces an involution $\psi^{[2]}$ on $S^{[2]}$. Furthermore the birational map $\psi^{[2]}\circ\varphi$ induces the identity on $H^2(S^{[2]},\mathbb{Z})$, therefore it is biregular (sends any K\"{a}hler class into itself), and it is also the identity (see \Ref{oss}{noidentity}). This means $\varphi=\psi^{[2]}$, which implies our claim.
\end{dim1}
\section*{Acknowledgements}
The author would like to thank his advisor, Prof. K.G. O'Grady, for the support in this work and for pointing out that the proofs of \Ref{lem}{denselemma} and \Ref{lem}{denselemma2} could be worked out from a proof on his notes on K3 surfaces. I would also like to thank Prof. E. Sernesi for his counsel concerning deformation theory, A. Rapagnetta for generous advice and Prof. S. Boissière for useful discussions.

\end{document}